\documentstyle[11pt, leqno, amssymb]{article}
    \normalbaselines\pretolerance=500	\tolerance=1000		\brokenpenalty=5000

\newcommand{\qed}{\sqcap\kern-8pt\sqcup}

\newenvironment{proof}
	{\par\noindent{\it Proof\/: }\nopagebreak\normalsize}%
	{\linebreak[2]\hspace*{\fill}$\qed$\ifdim\lastskip<12pt
	\removelastskip \penalty-200  \vskip12pt \fi}

\newtheorem{thm}{Theorem}[section]
\newtheorem{prop}[thm]{Proposition}
\newtheorem{lem}[thm]{Lemma}

\newtheorem{defn}[thm]{Definition}
\newtheorem{rem}[thm]{Remark}

\font\mathsymb=msbm10
\def\Bbb#1{\hbox{\mathsymb #1}}
\def\R{{\Bbb R}}

\def\ib{\begin{itemize}}
\def\ie{\end{itemize}}

\def \mod{{\;\mbox{\rm mod}}\;}

\def \rk{{\mbox{\rm rk}}}

\def \trace{{\mbox{\rm Trace}}\;}

\def \mod{{\;\mbox{\rm mod}}\;}

\def \rk{{\mbox{\rm rk}}}
\def \deg{{\mbox{\rm deg}\,}}

\def \sgn{{\mbox{\rm sgn}}}

\def \air{{\vskip 12pt\noindent}\par}

\def \ib{\begin{enumerate}}
\def \ie{\end{enumerate}}

\def \prs{{\mbox{\rm SRemS}}}

\def \sRes{{\mbox{\rm sRes}}}

\def \Id{{\mbox{\rm\bf Id}}}

\def \Ad{{\mbox{\rm\bf Ad}}}

\def \Bez{{\mbox{\rm Bez}}}

\def \TaQ{{\mbox{\rm TaQ}}}

\def \PmV{{\mbox{\rm PmV}}}

\def \Td{{\mbox{\rm Td}}}

\def \det{{\;\mbox{\rm det}}}
\def \sign{{\mbox{\rm sign}}}

\def \RR {{\Bbb R}}

\def\ib {\begin{enumerate}}
\def\ie {\end{enumerate}}

\def\adots{\mathinner{\mkern2mu\raise 1pt\hbox{.}\mkern 3mu\raise
4pt\hbox{.}\mkern1mu\raise 7pt\hbox{{.}}}}

\title{Sturm and Sylvester algorithms revisited via tridiagonal determinantal representations}
\author{Ronan Quarez\\ \small IRMAR (CNRS, URA 305), Universit\'e de Rennes 1, Campus de Beaulieu\\ 35042 Rennes Cedex, France\\ {\small e-mail : ronan.quarez@univ-rennes1.fr}}

\begin{document}
\maketitle

\begin{abstract}
\end{abstract}

{\small {\bf Keywords} : Determinantal representation - Hankel matrix - LU decomposition - Newton sums - Real roots - Signed remainders sequence - Sturm algorithm - Sylvester algorithm - Tridiagonal matrix}\par
{\small {\bf MSC Subject classification} : 12 - 15}
\begin{abstract} 
First, we show that Sturm algorithm and Sylvester algorithm, which compute the number of real roots of a given univariate polynomial, lead to two dual tridiagonal determinantal representations of the polynomial. Next, we show that the number of real roots of a polynomial given by a tridiagonal determinantal representation is greater than the signature of this representation.\end{abstract}

\section*{Introduction}
There are several methods to count the number of real roots of an univariate polynomial $p(x)\in\R[x]$ of degree $n$ (for details we refer to \cite{BPR}). Among them, the Sturm algorithm says that the number of real roots of $p(x)$ is equal to the number of Permanence minus the number of variations of signs which appears in the leading coefficients of the signed remainders sequence of $p(x)$ and $p'(x)$. \par
Another method is the Sylvester algorithm which says that the number of real roots of $p(x)$ is equal to the signature of the symmetric matrix whose $(i,j)$-th entry is the $i+j$-th Newton sums of the roots of the polynomial $p(x)$.\par 
One purpose of the paper is to point out, at least in the generic situation, that these two classical algorithms can be viewed as dual.\air

In section 1, we introduce signed remainders sequences of two given monic polynomials $p(x)$ and $q(x)$ of respective degrees $n$ and $n-1$. With some conventions of signs and others, we give a presentation of this sequence through a tridiagonal matrix $\Td(p,q)$. Next, we give a decomposition of this tridiagonal matrix as $\Td(p,q)=LC_p^{T}L^{-1}$ where $L$ is lower triangular and $C_p^T$ is the transposed of the companion matrix associated to $p(x)$. 
\par
In section 2, we introduce the duality between the Sturm and Sylvester algorithm, first when the polynomial $p(x)$ has only single and real roots, and then in Theorem \ref{duality} we generalize it to the generic case.\par
More precisely, on one hand we have $$\left\{\begin{array}{lcl}p(x)&=&\det(x\Id_n-\Td(p,q))\\ q(x)&=&\det(x\Id_{n-1}-\Td(p,q)_{n-1})\end{array}\right.$$  
with the conventions that $\Id_n$ (or $\Id$ in short) denotes the identity matrix of $\R^{n\times n}$ and $A_k\in\R^{k\times k}$ (respectively $\overline{A}_k\in\R^{k\times k}$) denotes the {\it $k$-th principal submatrix} (respectively {the {\it $k$-th antiprincipal submatrix}) of $A$ which corresponds to extracting the first $k$ (respectively the last $k$) rows and columns in the matrix $A\in\R^{n\times n}$. 

On the other hand, we consider a natural Hankel (hence symmetric) matrix $H(q/p)\in\R^{n\times n}$ associated to $p(x)$ and $q(x)$. Generically it admits an LU decomposition of the form $H(q/p)=KJK^T$ where $J$ is a {\it signature matrix} (a diagonal matrix with coefficients $\pm 1$ onto the diagonal) and $K$ is lower triangular. Then, we introduce the tridiagonal matrix $\overline{\Td}=K^{-1}C_p^TK$, which is such that $p(x)=\det(x\Id_n-\overline{\Td})$.\par
If we consider that the matrices $\Td(p,q)$ and $\overline{\Td}$ represent linear mappings in some basis, then the duality Theorem \ref{duality} means that one matrix can be deduced from the other simply by reversing the ordering of the basis.  

\air
We shall mention that, in the case when all the roots of $p(x)$ are real, the existence of a tridiagonal symmetric matrix $\Td$ given by the signed remainders sequence of $p(x)$ and $q(x)$ together with the identity $p(x)=\det(x\Id_n-\Td)$ corresponds to the Routh-Lanczos algorithm which answers a structured Jacobi inverse problem. Namely, the question to find a real symmetric tridiagonal matrix $A$ with a given characteristic polynomial $p(x)$ and such that the characteristic polynomial of its principal minor $A_{n-1}$, of size $n-1$, is proportional to $p'(x)$.  We refer to \cite{EP} for a survey on the subject.

\air
In section 3, we focus on the relation of the question of real roots counting and the question of determinantal representation. We say that $p(x)=\det(J-xA)$ is a determinantal representation of the polynomial $p(x)$ if $J\in\R^{n\times n}$ is a signature matrix and $A\in\R^{n\times n}$ is a symmetric matrix.\par
Remark that we may transform the identity $p(x)=\det(J-xA)$ into $p^*(x)=\det(xJ-A)$ where $p^*(x)$ is the reciprocal polynomial of $p(x)$. If we write $$p^*(x)=\det(J)\times\det(x\Id-AJ),$$ it shows a connexion with the results of section 2 when the matrix $AJ$ is tridiagonal. More precisely, we establish that such a determinantal representation is always possible : we may even find a family of representations for a given polynomial $p(x)$. We show also that, given such a determinantal representation for a polynomial $p(x)$, its number of real roots is at least equal to the signature of the signature matrix $J$.  \par

Finally, in section 4 we end with some worked examples.  

\section{Tridiagonal representation of signed remainders sequences}

\subsection{Definitions}\label{definitions}
Let $\alpha=(\alpha_1,\ldots,\alpha_n)$, $\beta=(\beta_1,\ldots,\beta_{n-1})$ and $\gamma=(\gamma_1,\ldots,\gamma_{n-1})$ be three sequences of real numbers. We set the tridiagonal matrix $\Td(\alpha,\beta,\gamma)$ to be :
 $$\Td(\alpha,\beta,\gamma)=\left(\begin{array}{ccccc}
\alpha_n&\gamma_{n-1}&0&\ldots&0\\
\beta_{n-1}&\alpha_{n-1}&\gamma_{n-2}&\ddots&\vdots\\
0&\beta_{n-2}&\ddots&\ddots&0\\
\vdots&\ddots&\ddots&\ddots&\gamma_1\\
0&\ldots&0&\beta_1&\alpha_1\\
\end{array}\right)
$$

Let $p(x)$ and $q(x)$ be two monic polynomials of respective degrees $n$ and $n-1$. We set $\prs(p,q)=(p_k(x))_k$ to be the {\it signed remainders sequence} of $p(x)$ and $q(x)$ defined in the following way :
\air
\begin{equation}\label{definition_prs}\left\{\begin{array}{lll}
p_0(x)&=&p(x)\\
p_1(x)&=&q(x)\\
p_k(x)&=&q_{k+1}(x)p_{k+1}(x)-\epsilon_{k+1}\beta^2_{k+1}p_{k+2}(x)
\end{array}\right.
\end{equation}
where 
\begin{equation}\label{convention_prs}\left\{\begin{array}{l}
p_k(x),q_{k+1}(x)\in\R[x],\\
\epsilon_{k+1}\in\{-1,+1\},\\
\beta_{k+1}\; {\rm is\; a\; positive\; real\; number},\\
p_{k+2}(x)\; {\rm is\; monic\; and\;} \deg p_{k+2}<\deg p_{k+1}.
\end{array}\right.
\end{equation}

This is a finite sequence which stops at the step just before we reach the zero polynomial as remainder. 
With these conventions, the signed remainders sequence $\prs(p,q)$ that we obtain is also called the {\it Sturm-Habicht} sequence of $p(x)$ and $q(x)$.\air

Let us assume that there is no degree breakdown in $\prs(p,q)$. Namely :
\begin{equation}\label{NoBreak}
(\forall k\in\{0,\ldots,n\})\; (\deg p_k=n-k)
\end{equation}
Then, $q_{k+1}(x)$ is a degree one polynomial which we write $q_{k+1}(x)=(x-\alpha_{k+1})$ with $\alpha_{k+1}\in\R$. Another consequence is that $\gcd(p,q)=1$.\air

Let $\gamma_{k+1}=\epsilon_{k+1}\beta_{k+1}$ and consider the following tridiagonal matrix :
$$\Td(p,q)=\Td(\alpha,\beta,\gamma)$$

We may read on this matrix all the informations about the signed remainders sequence $\prs(p,q)$. \air

For a given tridiagonal matrix $\Td=\Td(\alpha,\beta,\gamma)\in\R^{n\times n}$, we define the {\it first principal lower diagonal} (respectively the {\it first principal upper diagonal}) of $\Td$ to be the sequence $\beta=(\beta_1,\ldots,\beta_{n-1})$ (respectively $\gamma=(\gamma_1,\ldots,\gamma_{n-1})$). We will say that these first principal diagonals are {\it non-singular} if all the coefficients $\beta_i$ (respectively $\gamma_i$) are different from zero.\par Note that the no degree breakdown assumption (\ref{NoBreak}) implies that the principal diagonals of $\Td(p,q)$ are non-singular.  \par

\begin{prop}\label{prs_matrice_tridiag}
\ib
\item[(i)] To any tridiagonal matrix $\Td=\Td(\alpha,\beta,\gamma)$ with non-singular principal diagonals, we may canonically associate a (unique) couple of monic polynomials $p(x)$ and $q(x)$ of respective degrees $n$ and $n-1$ such that the sequence $\prs(p,q)$ has no degree breakdown and the characteristic polynomial of $\Td_k$ is equal to $p_{n-k}(x)$ : $$\det(x\Id_k-\Td_k)=p_{n-k}(x).$$ 

\item[(ii)] To any couple of monic polynomials $p(x)$ and $q(x)$ of respective degrees $n$ and $n-1$ such that $\prs(p,q)$ has no degree breakdown, we may associate a unique tridiagonal matrix with non-singular principal diagonals $\Td(p,q)=\Td(\alpha,\beta,\gamma)$ satisfying for all $k$, $\beta_k>0$ and $\gamma_k=\epsilon_k\beta_k$ where $\epsilon_k=\pm 1$.\par

 \item[(iii)] 
When we have $(i)$ and $(ii)$, the matrix $Td(p,q))\times P$ is tridiagonal and symmetric, where we have set $$P=\left(\begin{array}{rccc}
\epsilon_{n-1}\times\ldots\times\epsilon_1&&&\\
\ddots&&&\\
&\epsilon_2\times\epsilon_1&&\\
&&\epsilon_1&\\
&&&1
\end{array}\right).$$ 
 \item[(iv)] 
When we have $(i)$ and $(ii)$, the sequence of signs in the leading coefficients of the signed remainders sequence $\prs(p,q)$ is :
$$(1,1,\epsilon_1,\epsilon_2,\epsilon_1\times\epsilon_3,\epsilon_2\times\epsilon_4,\epsilon_1\times\epsilon_3\times\epsilon_5,\ldots,\epsilon_{n-1\mod 2}\times\ldots\times\epsilon_{n-3}\times\epsilon_{n-1})$$
\ie
\end{prop} 
 
\begin{proof}
Concerning $(i)$, the polynomials $p(x)$ and $q(x)$ are taken to be $p(x)=\det(x\Id_n-\Td)$ and $q(x)=\det(x\Id_{n-1}-\Td_{n-1})$. Then, we set for all $k$, $$\delta_{n-k}(x)=\det(x\Id_k-\Td_k)$$ (where $\Td_k$ is the $k$-th principal submatrix of $\Td$) and we develop the determinant $$\delta_0(x)=\det(x\Id_n-\Td)$$ with respect to the last row. We get 
$$\delta_0(x)=(x-\alpha_1)\delta_1(x)-(\beta_1\gamma_1)\delta_2(x)$$
Repeating the process, we obtain the same recurrence relation as the one defining the sequence $(p_k(x))_k$ in (\ref{definition_prs}). Since $\delta_{0}(x)=p_{0}(x)$ and $\delta_{1}(x)=p_{1}(x)$, we get the wanted identity.\air
Point $(ii)$ follows straightforward from the beginning of the section, whereas points $(iii)$ and $(iv)$ follows from elementary computation. \end{proof}

We may note that to the tridiagonal matrix $\Td(p,q)$, we may associate also another natural polynomial remainder sequence : $\overline{\prs(p,q)}=\prs(p,\bar{q})$ where $$p(x)=\det(x\Id_n-\Td)$$ and $$\bar{q}(x)=\det(x\Id_{n-1}-\overline{\Td}_{n-1}),$$ with the convention that $\overline{\Td}_k$ is the $k$-th antiprincipal submatrix of $\Td$.
The signed remainders sequence $\overline{\prs(p,q)}$ will be considered as the {\it dual} signed remainders sequence of $\prs(p,q)$.
This only means that we may read on a tridiagonal matrix from the top left rather than from the bottom right !

\air
For cosmetic reasons we will write $\overline{\Td(p,q)}$ in place of $\Td(p,\bar{q})$. We obviously have :

\begin{equation}\label{dual}\overline{\Td(p,q)}=\Ad\times\Td(p,q)\times\Ad\end{equation}

where $\Ad_n\in\R^{n\times n}$ ($\Ad$ in short) stand for the anti-identity matrix of size $n$ :

$$\Ad_n=\left(\begin{array}{cccc}
0&\ldots&0&1\\
\vdots&\adots&\adots&0\\
0&\adots&\adots&\vdots\\
1&0&\ldots&0\\
\end{array}\right)$$ 

\subsection{Companion matrix}\label{Companion}
We denote by $A^T$ the transposed of the matrix $A\in\R^{n\times n}$ and we define the {\it companion matrix} of the polynomial $p(x)=x^n+a_{n-1}x^{n-1}+\ldots+a_0$ to be
$$C_p=
\left(\begin{array}{ccccc}
0&\ldots&\ldots&0&-a_0\\
1&\ddots&&\vdots&-a_1\\
0&\ddots&\ddots&\vdots&\vdots\\
\vdots&\ddots&\ddots&0&-a_{n-2}\\
0&\ldots&0&1&-a_{n-1}
\end{array}\right)
$$

We recall a well-know identity (see for instance \cite{EP}) :
\begin{prop}\label{prs_tridiag} Let $p(x)$ and $q(x)$ be two monic polynomials of respective degrees $n$ and $n-1$ such that $\prs(p,q)$ has no degree breakdown. \par 
Then, there is a lower triangular matrix $L$ such that \begin{equation}\Td(p,q)=LC_p^TL^{-1}\end{equation}
\end{prop}

\begin{proof}
With the notation of Subsection \ref{definitions}, let ${\cal P}(x)=\left(\gamma_1\ldots\gamma_{n-1}p_n(x),\ldots,\gamma_1p_2(x),p_1(x)\right)$.
A direct computation gives
$${\cal P}(x)\left(\Td(p,q)\right)^T=x{\cal P}(x)+\left(0,\ldots,0,-p(x)\right)$$
Let $U$ be the upper triangular matrix whose columns are the coefficients of the polynomials of ${\cal P}(x)$ in the canonical basis ${\cal C}(x)=(1,x,\ldots,x^{n-1})$.
In other words :
$${\cal C}(x)U={\cal P}(x)$$

Besides, we have $${\cal C}(x)C_p=x{\cal C}(x)+(0,\ldots,0,-p(x))$$
Thus
$$\begin{array}{rclcc}
{\cal C}(x)C_pU&=&x{\cal C}(x)U+(0,\ldots,0,-p(x))&{\rm since}\; p_1(x) \;{\rm is\; monic}\\
&=&{\cal P}(x)\left(\Td(p,q)\right)^T&\\
&=&{\cal C}(x)U\left(\Td(p,q)\right)^T&\\
\end{array}$$

We deduce the identity 
$$V(x_1,\ldots,x_n)C_pU=V(x_1,\ldots,x_n)U\left(Td(p,q)\right)^T$$
for any Vandermonde matrix $V(x_1,\ldots,x_n)$ whose lines are $(1,x_i,\ldots,x_i^{n-1})$ for $i=1\ldots n$. If we choose the $n$ reals $x_1,\ldots,x_n$ to be distinct, then $V(x_1,\ldots,x_n)$ becomes invertible and we get :
$$\Td(p,q)=LC_p^TL^{-1}$$
where $L$ is the lower triangular matrix defined by $L= U^T$.
\end{proof}

The following result says that the decomposition generically exists for any tridiagonal matrix, and also it is unique :
\begin{prop}\label{commutant}
Any tridiagonal matrix $\Td$ with non-singular principal diagonals can be written $\Td=LC_p^TL^{-1}$ where  $p(x)=\det(x\Id-\Td)$ and $L$ is a lower triangular matrix. Moreover the matrix $L$ is unique up to a multiplication by a real number.  
\end{prop}
\begin{proof}
The existence is given by Proposition \ref{prs_matrice_tridiag} and Proposition \ref{prs_tridiag}.\par
We come now to the unicity. 
Assume that $L_1C_p^TL_1^{-1}=L_2C_p^TL_2^{-1}$ where $L_1$ and $L_2$ are lower triangular. Then, $L=L_2^{-1}L_1$ is a lower triangular matrix which commute with $C_p^T$. 
\par
If $L=(t_{i,j})_{1\leq i,j\leq n}$, then 
$$LC_p^T=\left(\begin{array}{ccccc}
0&t_{1,1}&0&\ldots&0\\
\vdots&t_{2,1}&t_{2,2}&\ddots&\vdots\\
\vdots&\vdots&&\ddots&0\\
0&t_{n-1,1}&\ldots&\ldots&t_{n-1,n-1}\\
?&\ldots&\ldots&\ldots&?
\end{array}\right)$$
and
$$C_p^TL=\left(\begin{array}{ccccc}
t_{2,1}&t_{2,2}&0&\ldots&0\\
t_{3,1}&t_{3,2}&t_{3,3}&\ddots&\vdots\\
\vdots&&\ddots&\ddots&0\\
t_{n,1}&\ldots&\ldots&t_{n,n-1}&t_{n,n}\\
?&\ldots&\ldots&\ldots&?
\end{array}\right)$$
\air
Thus $t_{1,1}=t_{2,2}=\ldots=t_{n,n}$ and $t_{2,1}=t_{3,2}=\ldots=t_{n,n-1}=0$ and $t_{3,1}=t_{4,2}=\ldots=t_{n,n-2}=0$, and so on until $t_{n,1}=0$. We deduce that $L=\lambda \Id$ and we are done. 
\end{proof}

\subsection{Sturm algorithm}
As a particularly important case of signed remainders sequences, we shall mention the Sturm sequence which is $\prs(p,q)$ where $q$ is taken to be the derivative of the polynomial $p(x)$ up to normalization, i.e. $q=p'/\deg(p)$. \par

For a given finite sequence $\nu=(\nu_1,\ldots,\nu_k)$  of elements in $\{-1,+1\}$, we recall the {\it Permanence minus Variations} number :
$$\PmV(\nu_1,\ldots,\nu_k)=\sum_{i=1}^{k-1}\nu_i\nu_{i+1}.$$
Here the sequence $\nu$ will be for the sequence of signs of leading coefficients in $\prs(p,q)$. Then, the Sturm Theorem \cite[Theorem 2.50]{BPR} says that the number $\PmV(\nu)$ is exactly the number of real roots of $p(x)$. \par

If we assume that the polynomial $p(x)$ has $n$ distinct real roots, then the Sturm sequence has no degree breakdown and for all $k$ we have $\nu_k=1$. Hence we get a {\it symmetric} tridiagonal matrix $\Td(p,q)$ which has the decomposition 
$\Td(p,q)=LC_p^TL^{-1}$ where $L$ is the lower triangular matrix defined as in subsection \ref{Companion}. In particular, the last row of $L$ gives the list of coefficients of the polynomial $q(x)$ in the canonical basis.

\section{Duality between Sturm and Sylvester algorithms}

\subsection{Sylvester algorithm}\label{Sylvester}

Let us introduce the symmetric matrix ${\rm Newt}_p(n)=(n_{i,j})_{0\leq i,j\leq n-1}$ define as $$n_{i,j}=\trace(C_p^{i+j})=N_{i+j}$$ which is nothing but the $i+j$-th Newton sum of the polynomial $p(x)$. To be more explicit, if $\alpha_1,\ldots,\alpha_n$ denote all the complex roots of the polynomial $p(x)$, then the $k$-th Newton sum is the real number $N_k=\alpha_1^k+\ldots+\alpha_n^k$.  
\air
Recall that the {\it signature} $\sign(A)$ of a real symmetric matrix $A\in\R^{n\times n}$, is defined to be the number $p-q$, where $p$ is the number of positive eigenvalues of $A$ (counted with multiplicity) and $q$ the number of negative eigenvalues of $A$ (counted with multiplicity).
The Sylvester Theorem (which has been generalized later by Hermite :  \cite[Theorem 4.57]{BPR}) says that the matrix ${\rm Newt}_p(n)$ is invertible if and only if $p(x)$ has only single roots, and also that $\sgn({\rm Newt}_p(n))$ is exactly the number of distinct real roots of $p(x)$. \par
In particular, if the polynomial $p(x)$ has $n$ distinct real roots, then the matrix ${\rm Newt}_p(n)$ is positive definite. 
Thus, by the Choleski decomposition algorithm, we can find a lower triangular matrix $K$ such that ${\rm Newt}_p(n)=K K^T$. Le us show how to exploit this decomposition.\air
First, we write $$p(x)=\det(x\Id-C_p^T)$$ 
Then, we introduce a useful identity (which will be discussed in more details in the forthcoming section) : 
$${\rm Newt}_p(n)C_p=C_p^T{\rm Newt}_p(n),$$ 
So, we get :
$$p(x)=\det(x\Id-K^{-1}C_p^TK)$$ 
Note that the matrix $K^{-1}C_p^TK$ is tridiagonal. Our purpose in the following is to establish a connexion with the identity $$p(x)=\det(x\Id-LC_p^TL^{-1})$$
obtained in Proposition \ref{commutant}.\par
More generally, we will point out a connexion between tridiagonal representations associated to signed remainders sequences on one hand, and tridiagonal representations derived from decompositions of some Hankel matrices on the other hand.

\subsection{Hankel matrices and Intertwinning relation}
Roughly speaking, the idea of previous section is to start with the canonical companion identity $$p(x)=\det(x\Id-C_p^T)$$ and then to use a symmetric invertible matrix $H$ satisfying the so-called {\it intertwinning relation}

\begin{equation}\label{Intertwinning} HC_p=C_p^TH\end{equation}

Since $H$ is supposed to be symmetric invertible, Equation (\ref{Intertwinning}) only says that the matrix $HC_p$ is symmetric. It is a classical and elementary result that a matrix $H$ satisfying equation (\ref{Intertwinning}) is necessarily an Hankel matrix. 

\begin{defn}
We say that the matrix $H=(h_{i,j})_{0\leq i,j\leq n-1}\in\R^{n\times n}$ is an Hankel matrix if $h_{i,j}=h_{i',j'}$ whenever $i+j=i'+j'$. Then, it makes sense to introduce the real numbers $a_{i+j}=h_{i,j}$ which allow to write in short $H=(a_{i+j})_{0\leq i,j\leq n-1}$.
\air
Let $s=(s_k)$ be a sequence of real numbers. We denote by $H_n(s)$ or by $H(s_0,\ldots,s_{2n-2})$ the following Hankel matrix of $\R^{n\times n}$ :

$$H_n(s)=(s_{i+j})_{0\leq i,j\leq n-1}=\left(\begin{array}{cccc}
s_0&s_1&\ldots&s_{n}\\
s_1&&\adots&s_{n+1}\\
\vdots&\adots&\adots&\vdots\\
s_n&s_{n+1}&\ldots&s_{2n-2}\\
\end{array}\right)$$
\air
\end{defn}

We get from \cite[Theorem 9.17]{BPR} :
\begin{prop}\label{Hankel} Let $p(x)=x^n+a_{n-1}x^{n-1}+\ldots+a_0$ and $s=(s_k)$ be a sequence of real numbers. The following assertions are equivalent
\ib
\item[(i)]
$(\forall k\geq n)\; (s_k=-a_{n-1}s_{k-1}-\ldots-a_0s_{k-n})$

\item[(ii)]
There is a polynomial $q(x)$ of degree $\deg q<\deg p$ such that $$\frac{q(x)}{p(x)}=\sum_{j=0}^\infty\frac{s_j}{x^{j+1}}$$
\item[(iii)]
There is an integer $r\leq n$ such that $\det(H_r(s))\not =0$, and for all $k>r$, $\det(H_k(s))=0$.
\ie
Whenever these conditions are fulfilled, we denote by $H_n(q/p)$ the Hankel matrix $H_n(s)$.
\end{prop}

Back to the intertwinning relation (\ref{Intertwinning}) : it is immediate that an Hankel matrix $H$ is a solution if and only if the (finite) sequence $(s_0,\ldots,s_{2n-2})$ satisfies the linear recurrence relation of Proposition \ref{Hankel}$(i)$, for $k=n,\ldots,2n-2$.\par 
For further details and developments about the intertwinning relation,  
we refer to \cite{HV}.\air

The vector subspace of Hankel matrices in $\R^{n\times n}$ satisfying relation (\ref{Intertwinning}) has dimension $n$, and contains a remarkable element : the Hankel matrix ${\rm Newt}_p(n)$ that was considered in subsection \ref{Sylvester} about Sylvester algorithm.
Indeed, it is a well-known and elementary fact that the $N_k$'s are real numbers which verify the Newton identities :
$$(\forall k\geq n)\;\left(N_k+a_{n-1}N_{k-1}+\ldots+a_0N_{k-n}=0\right)$$    
\air

\subsection{Barnett formula}
First, recall that if $p(x)=x^n+a_{n-1}x^{n-1}+\ldots+a_0$ and $q(x)$ is  a (non-necessarily monic) polynomial in $\R[x]$ whose degree is equal to $n-1$, the Bezoutian of $p(x)$ and $q(x)$ is defined as the two-variables polynomial :

$$\Bez (p,q)=\frac{q(y)p(x)-q(x)p(y)}{x-y}\in\R[x,y]$$

Let ${\cal B}(z)$ be any basis of the $n$-dimentional vector space $\R[z]/p(z)$ over $\R$. We denote by $\Bez_{\cal B}(p,q)$ the symmetric matrix of the coefficients of $\Bez(p,q)$ 
with respect to the basis ${\cal B}(x)$ and  ${\cal B}(y)$.\par
Among all the basis of $\R[z]/p(z)$ that will be interesting for the following, let us mention the canonical basis ${\cal C}=(1,z,\ldots,z^{n-1})$ and also the (degree decreasing) Horner basis ${\cal H}(z)=(h_0,\ldots,h_{n-1})$ associated to the polynomial $p(z)$ and which is defined by : 
$$\left\{\begin{array}{lcl}
h_{0}(z)&=&z^{n-1}+a_{n-1}z^{n-2}+\ldots+a_{1}\\
&\vdots&\\
h_{i}(z)&=&z^{n-1-i}+a_{n-1}z^{n-2-i}+\ldots+a_{i+1}=zh_{i+1}(z)+a_{i+1}\\
&\vdots&\\
h_{n-2}(z)&=&z+a_{n-1}\\
h_{n-1}(z)&=&1
\end{array}\right.$$

\air
We recall from \cite[Proposition 9.20]{BPR} :
\begin{prop}\label{BPR9.20} Let $p(x)$ and $q(x)$ be two polynomials such that $\deg q<\deg p=n$. Let $s$ be the sequence of real numbers defined by $$\frac{q(x)}{p(x)}=\sum_{j=0}^\infty\frac{s_j}{x^{j+1}}$$
Then, $\Bez_{\cal H}(p,q)=H_n(s)=H_n(q/p)$.
\end{prop}
\air

\air We come to a central proposition which is a consequence of the Barnett formula \cite{Ba}. 

\begin{prop}\label{Barnett} Let $p(x)$ and $q(x)$ be two polynomials such that $\deg q<\deg p=n$ and let $P_{\cal C H}$ be the change of basis matrix from the canonical basis ${\cal C}$ to the Horner basis ${\cal H}$. We have 
$$q(C_p)=P_{{\cal CH}}^T\times H_n(q/p)$$
\end{prop}

\begin{proof}
The Barnett formula has been established in \cite{Ba} using direct matrix computations. For the convenience of the reader, we give here another proof (which may be found at various places in the literature). 

The obvious identity 
$$q(y)(p(x)-p(y))=q(y)p(x)-p(y)q(x)+p(y)(q(x)-q(y))$$
implies, by definition of the Bezoutian $B(p,q)$, that : 

$$q(y)\frac{p(x)-p(y)}{x-y}\equiv \Bez(p,q)\mod p(y)$$
Noticing that $\frac{p(x)-p(y)}{x-y}=\sum_{j=0}^{n-1}h_j(y)x^j$, we get 
$$q(y)\sum_{j=0}^{n-1}h_j(y)x^j\equiv {\cal C}(y) \Bez_{\cal C}(p,q){\cal C}(x)^T\mod p(y)$$

In other words, if we denote by $M$ the matrix whose columns are the coefficients of $q(y)h_j(y)$ in the basis ${\cal C}(y)$, we get the identity 
$${\cal C}(y) M {\cal C}(x)^T \equiv {\cal C}(y) \Bez_{\cal C}(p,q){\cal C}(x)^T\mod p(y)$$

Since $C_p$ is the matrix of the multiplication by $y\mod p(y)$ with respect to the canonical basis ${\cal C}(y)$, we have also the identity
$$M=q(C_p)P_{\cal C H}$$
where the change of basis matrix $P_{\cal C H}$ is in fact the following Hankel matrix $$P_{\cal C H}=H(a_0,a_1,\ldots,a_{n-1},1,0,\ldots,0)\in\R^{n\times n}$$
with the usual notation $p(x)=x^n+a_{n-1}x^{n-1}+\ldots+a_0$. 
Hence, we get the Barnett Formula

$$\Bez_{\cal C}(p,q)=q(C_p)P_{\cal C H}$$

Finally, by Proposition \ref{BPR9.20}, we derive the wanted relation :
$$\Bez_{\cal C}(p,q)=P_{\cal C H}^T\Bez_{\cal H}(p,q)P_{\cal C H}=P_{\cal C H}^TH_n(q/p)P_{\cal C H}$$

Which concludes the proof.
\end{proof}

To end the section, we show how Sturm and Sylvester algorithms can be considered as dual, in the case where all the roots of $p(x)$ are real and simple, say $x_1<\ldots<x_n$.  Then, $q(x)=p'(x)/n$ has also $n-1$ simple real roots $y_1<\ldots<y_{n-1}$ which are interlacing those of $p(x)$. Namely 
$$x_1<y_1<x_2<y_2<\ldots<y_{n-1}<x_n$$ 
We may repeat the argument to see that this interlacing property of real roots remains for any two consecutive polynomials $p_k(x)$ and $p_{k+1}(x)$ of the sequence $\prs(p,q)$. 
In particular, $\prs(p,q)$ does not have any degree breakdown, all the $\epsilon_k$ are equal to $+1$, and $N(q/p)$ is positive definite. \air

We have, by Proposition \ref{Barnett} $$q(C_p^T)=H_n(q/p)P_{\cal CH}$$ 
Since $H_n(q/p)$ is positive definite, the Cholesky algorithm gives a decomposition $$H_n(q/p)= K K^T$$ where $K\in\R^{n\times n}$ is lower triangular. So that we can write
$$p(x)=\det(x\Id-K^{-1}C_p^TK)$$ 
We shall remark at this point that the matrix $K^{-1}C_p^TK$ is tridiagonal and symmetric.

\air

We get 
$q(C_p^T)=K\Ad L$ where $L=\Ad K^TP_{\cal CH}$.
Then, we observe that $L$ is a lower triangular matrix (since $P_{\cal CH}\Ad$ is upper triangular) and $K\Ad L$ commute with $C_p^T$.
Thus, we have the identity : 
$$LC_p^TL^{-1}=\Ad(K^{-1}C_p^TK)\Ad$$
We denote by $\Td$ this tridiagonal matrix. Let $(p_k(x))$ be the signed remainders sequence associated to $\Td$ as given in Proposition \ref{prs_matrice_tridiag} (i). The first row of $K\Ad L$ is proportional to the last row of the matrix $L$ which is proportional to $p_{1}(x)$.  
It remains to observe that the first row of $K\Ad L=q(C_p^T)$ gives exactly the coefficients of the polynomial $q(x)$ in the canonical basis. 
Then, $p_1(x)=q(x)$. \air

In summary, we have shown that, if $p(x)$ has $n$ simple real roots and $q(x)=p'(x)/n$, then $H_n(q/p)$ is positive definite with Cholesky decomposition $H_n(q/p)=KK^T$, and if we denote by $\tilde{q}(x)$ the monic polynomial whose coefficients are proportional to the last row of $K^{-1}$, then $\Td(p,\tilde{q})=\overline{\Td(p,q)}$. Which settle the announced duality.

\subsection{Generic case}
We turn now to the generic situation. Let $p(x)$ and $q(x)$ be monic polynomials of respective degrees $n$ and $n-1$ and such that $\prs(p,q)$ does not have any degree breakdown. This condition is equivalent to saying that all the principal minors of the Hankel matrix $H_n(q/p)$ do not vanish. We refer to \cite{BPR} for this point. One way to see this is to figure out the connexion with the subresultants of $p(x)$ and $q(x)$.\par A little bit more precisely, the $j$-th signed subresultant coefficient of $p(x)$ and $q(x)$ is denoted by $\sRes_j(p,q)$ for $j=0\ldots n-1$. If for all $j$, $\sRes_j(p,q)\not = 0$, we say that the sequence of subresultants is {\it non-defective}.  
Then, by \cite[Corollary 8.33]{BPR} and Proposition \ref{prs_matrice_tridiag} (iv), we deduce that the non-defective condition is equivalent to the fact that $\prs(p,q)$ has no degree breakdown. Moreover, from \cite[Lemma 9.26]{BPR} we know that $$(\forall j\in\{1\ldots n\})\;\left(\sRes_{n-j}(p,q)=\det(H_j(q/p))\right).$$
In conclusion, our no degree breakdown assumption means also that all the principal minors of the Hankel matrix $H_n(q/p)$ do not vanish.

\par At this point, we may add another equivalent condition, which will be essential for the following. Indeed, the condition that all the principal minors of the Hankel matrix $H_n(q/p)$ do not vanish is also equivalent to saying that the matrix $H_n(q/p)$ admits an invertible $LU$ decomposition. Namely, it exists a lower triangular matrix $L$ with $1$ entries onto the diagonal, and an upper invertible triangular matrix $U$ such that  $H_n(q/p)=LU$. Moreover this decomposition is unique and since $H_n(q/p)$ is symmetric we may write it as $H_n(q/p)=LDL^T$ where $D$ is diagonal. In fact, for our purpose, we will prefer the unique decomposition $H_n(q/p)=KJK^T$ where $K$ is lower triangular and $J$ is a signature matrix. 

Generalizing the previous section, we get :

\begin{thm}\label{duality}
Let $p(x)$ and $q(x)$ be two monic polynomials of respective degrees $n$ and $n-1$ such that $\prs(p,q)$ does not have any degree breakdown. Consider the symmetric LU-decomposition of the Hankel matrix $H_n(q/p)=KJK^T$, where $J$ is a signature matrix and $K$ a lower triangular matrix, and denote by $\tilde{q}(x)$ the monic polynomials whose coefficients in the canonical basis are proportional to the last row of $K^{-1}$. Then, 
$$\Td(p,\tilde{q})=\overline{\Td(p,q)}.$$
\end{thm}

\begin{proof}
We start with the companion identity :
$$p(x)=\det(x\Id-C_p^T)$$
Next, because of Proposition \ref{Hankel}$(i)$, we notice that the matrix $H_n(q/p)$ verifies the intertwinning relation :
  
$$H_n(q/p)C_p=C_p^TH_n(q/p)$$

Then, we write the symmetric $LU$-decomposition of $H_n(q/p)$ :
$$H_n(q/p)=KJK^T$$
Which gives the identity 
 
$$p(x)=\det(x\Id-K^{-1}C_p^TK).$$ 
\air

We have, by Proposition \ref{Barnett} $$q(C_p^T)=H_n(q/p)P_{\cal CH}=K\Ad L$$ 
where $$L=\Ad J K^TP_{\cal CH}$$

We observe first that $L$ is a lower triangular matrix (since $P_{\cal CH}\Ad$ is upper triangular), and second that $K\Ad L$ commute with $C_p^T$.
Thus, we have the identity : 
$$LC_p^TL^{-1}=\Ad(K^{-1}C_p^TK)\Ad$$

Proposition \ref{commutant} gives 
$$\overline{\Td(p,\tilde{q})}=\Ad(K^{-1}C_p^TK)\Ad$$

Moreover, the first row of $K\Ad L$ is proportional to the last row of the matrix $L$. It remains to observe that the first row of $K\Ad L=q(C_p^T)$ gives exactly the coefficients of the polynomial $q(x)$ in the canonical basis. Thus, by Proposition \ref{commutant} we get

$$LC_p^TL^{-1}=\Td(p,q)$$
Which concludes the proof.
\end{proof}

\begin{rem}
Note that $K^{-1}C_p^TKJ$ is symmetric and hence also the matrix $LC_p^TL^{-1}\bar{J}$, where $\bar{J}=\Ad J\Ad$.    
\end{rem}

\section{Tridiagonal determinantal representations}

\subsection{Notations}
We say that an univariate polynomial $p(x)\in\R[x]$ of degree $n$ such that $p(0)\not =0$ has a {\it determinantal representation} if 
   
$${\rm (DR)}\quad p(x)=\alpha\det(J-Ax)$$ 

where $\alpha\in\RR^*$, $J$ is a signature matrix in $\RR^{n\times n}$, and $A$ is a symmetric matrix in $\RR^{n\times n}$ (we obviously have $\alpha=\det(J)p(0)$). \par

Likewise, we say that $p(x)$ has a {\it weak determinantal representation} if
   
$${\rm (WDR)}\quad p(x)=\alpha\det(S-Ax)$$ 

where $\alpha\in\R^*$, $S$ is symmetric invertible and $A$ is symmetric.\air

Of course the existence of (DR) is obvious for univariate polynomials, but we will focus on the problem of {\it effectivity}. Namely, we want an algorithm (say of polynomial complexity with respect to the coefficients and the degree of $p(x)$) which produces the representation. Typically, we do want to avoid the use of the roots of $p(x)$.

One result in that direction can be found in \cite{Qz2} (which is inspired from \cite{Fi}). It uses arrow matrices as a ``model", whereas in the present article we make use of tridiagonal matrices.\air 

When all the roots of $p(x)$ are real, the effective existence of determinantal representation for univariate real polynomials exists even if we add the condition that $J=\Id$. It has been discussed in several places, although not exactly with the determinantal representation formulation. Indeed, in place of looking for DR we may consider the equivalent problem of the research of a symmetric matrix whose characteristic polynomial is given. Indeed, if the size of the matrix $A$ is equal to the degree $n$ of the polynomial, the condition $$p(x)=\det(\Id-xA)$$ is equivalent to $$p^*(x)=\det(x\Id-A)$$ where $p^*(x)$ is the reciprocal polynomial of $p(x)$.
In \cite{Fi}, arrow matrices are used  to answer this last problem. On the other hand, the Routh-Lanczos algorithm (which can be viewed as Proposition \ref{prs_matrice_tridiag}) gives also an answer, using  tridiagonal model. Note that the problem may also be reformulated as a structured Jacobi inverse problem (confer \cite{EP} for a survey). \air
In the following, we generalize the tridiagonal model to any polynomial $p(x)$, possibly having non real roots. Doing that, general signature matrices $J$ appear, whose entries depend on the number of real roots of $p(x)$.\par

\subsection{Over a general field}
A lot of identities in Section 2 are still valid over a general field $k$.
For instance, if $p(x)$ and $q(x)$ are monic polynomials of respective degrees $n-1$ and $n$, we may still associate the Hankel matrix $H(q/p)=(s_{i+j})_{0\leq i,j\leq n-1}\in k^{n\times n}$ defined by the identity
$$\frac{q(x)}{p(x)}=\sum_{j=0}^\infty\frac{s_j}{x^{j+1}}$$ 

Then, we have the following :

\begin{thm}\label{general_field}
Let $p(x)\in k[x]$, $q(x)\in k[x]$ be two monic polynomials of respective degrees $n$ and $n-1$, and set $H=H(q/p)$. Then, the matrix $C_p^TH$ is symmetric and we have the WDR :  
$$\det(H)\times p(x)=\det(xH-C_p^TH)$$ 
Moreover, if $H$ admits the LU-decomposition $H=KDK^T$ where $K\in k^{n\times n}$ is lower triangular with entries $1$ onto the diagonal and $D\in k^{n\times n}$ a diagonal matrix, then we have : 
$$p(x)=\det(xD-\Td)$$ 
where $\Td=K^{-1}C_p^TKD$ is a tridiagonal symmetric matrix.
\end{thm}

\begin{proof}
We exactly follow the proof of Theorem \ref{duality}.
\end{proof}

Note that the condition for $H$ to be invertible is equivalent to the fact that the polynomials $p(x)$ and $q(x)$ are coprime, since we have 
$$\rk(\Bez(q,p))=\deg(p)-\deg(\gcd(p,q)).$$ 
To see this, we may refer to the first assertion of \cite[Theorem 9.4]{BPR} whose proof is valid over any field.\air

The WDR of Theorem \ref{general_field} has the advantage that the considered matrices have entries in the ring generated by the coefficients of the polynomial $p(x)$. 
This point is not satisfied in the methods proposed in $\cite{Qz2}$ or in the Routh-Lanczos algorithm.\par 
In fact, the use of Hankel matrices satisfying the intertwinning relation seems to be more convenient since we are able to ``stop the algorithm at an earlier stage", namely before having to compute a square root of the matrix $H$ (or of the matrix $D$).\par 

Of course, at the time we want to derive a DR, then we have to add some conditions on the field $k$, for instance we shall work over an ordered field where square roots of positive elements exist.\par 
To end the section, we may summarize that, for a given polynomial $p(x)$, we have an obvious but non effective (i.e. using factorization) DR with entries in the splitting field of $p(x)$ over $k$, to compare with
an effective WDR given by Theorem \ref{general_field} where entries are in the field generated by the coefficients of $p(x)$.    

\subsection{Symmetric tridiagonal representation and real roots counting}
If $p(x)$ and $r(x)$ are two real polynomials, we recall the number known as the Tarski Query :
$$\TaQ(r,p)=\#\{x\in\R\mid p(x)=0\wedge r(x)>0\}-\#\{x\in\R\mid p(x)=0\wedge r(x)<0\}.$$
We also recall the definition of the Permanences minus variations number of a given sequence of signs $\nu=(\nu_1,\ldots,\nu_k)$ :
$$\PmV(\nu)=\sum_{i=1}^{k-1}\nu_i\nu_{i+1}.$$

We summarize, from \cite[Theorem 4.32, Proposition 9.25, Corollary 9.8]{BPR} some useful properties of these numbers, 
\begin{prop}\label{real_counting}
Let $p(x)$ and $q(x)$ be two monic polynomials of respective degrees $n$ and $n-1$, and such that the sequence $\prs(p,q)$ has no degree breakdown. Let $r(x)$ be another polynomial such that $q(x)$ is the remainder of $p'(x)r(x)$ modulo $p(x)$. Then, 
$$\PmV(\nu)=\sgn(\Bez(p,q))=\sgn(H_n(q/p))=\TaQ(r,p)$$
where $\nu$ is the sequence of signs of the leading coefficients in the signed remainders sequence $\prs(p,q)$. 
\end{prop}

We come now to our main result about real roots counting :

\begin{thm}\label{realroots_tridiag}
Let $\Td\in\R^{n\times n}$ be a tridiagonal symmetric matrix with non-singular first principal diagonals. Let also $p(x)\in\R[x]$ be a real polynomial with no multiple root and such that $$p(x)=\det(J)\det(xJ-\Td),$$ where $J$ is a signature matrix whose last entry onto the diagonal is $+1$.\par 
Then, the number of real roots of $p(x)$ is greater than $\sgn(J)$. 
\end{thm}
\begin{proof} 
We have $$p(x)=\det(x\Id_n-\Td\times J)$$ and we set 
 $$q(x)=\det\left(x\Id_{n-1}-\left(\Td\times J\right)_{n-1}\right).$$
The matrix $\Td\times J$ is still tridiagonal with non-singular first principal diagonals. Then, we consider the sequence $\prs(p,q)$ and denote by $\nu$ the associated sequence of signs of leading coefficients.\par
Since $\gcd(p,p')=1$, we set $r(x)$ to be the unique polynomial of degree $<n$ such that $$r\equiv \frac{q}{p'}\mod p.$$ Then, $$p'r\equiv q\mod p$$ and from Proposition \ref{real_counting}, we get :
$$\PmV(\nu)=\TaQ(r,p)\leq \#\{x\in\R\mid p(x)=0\}$$

\air
Let us introduce some notations at this step. First, let $\Td=\Td(\alpha,\beta,\gamma)$, next denote by $\epsilon(a)$ the sign in $\{-1,+1\}$ of the non zero real number $a$, and finally let 
$$J=\left(\begin{array}{cccc}
\theta_{n-1}&&&\\
&\ddots&&\\
&&\theta_{1}&\\
&&&1
\end{array}\right).$$ 

Then, we can write

$$p(x)=\det\left(x\Id_n-P(\Td\times J)P^{-1}\right)$$ where 

$$P=\left(\begin{array}{ccccc}
(\theta_{n-1}\ldots\theta_{1})\times(\epsilon(\gamma_{n-1})\ldots\epsilon(\gamma_{1}))&&&&\\
&&\ddots&&\\
&&&\theta_{1}\times\epsilon({\gamma_1})&\\
&&&&1
\end{array}\right).$$

We note in fact that $P(\Td\times J)P^{-1}=\Td(p,q)$. Indeed, all the coefficients onto the first lower principal diagonal are positive. Moreover, all the coefficients onto the first upper principal diagonal are given by the sequence
 $$(\theta_{n-1}\times\theta_{n-2},\ldots,\theta_{2}\times\theta_1,\theta_1).$$
We deduce from Proposition \ref{prs_matrice_tridiag} (iv) that the sequence of signs of leading coefficients in the signed remainders sequence $\prs(p,q)$ is the following :
 $$\nu=(\theta_{n-1}\times\ldots\times\theta_{1},\ldots,\theta_2\times\theta_1,\theta_{1},1,1).$$
Thus $$\PmV(\nu)=1+\sum_{k=1}^{n-1}\theta_k=\sgn(J)$$ and we are done.
\end{proof}

Another way, maybe less constructive, to prove the result is to use the duality of Theorem \ref{duality}. Indeed, replacing as in the previous proof the matrix $\Td\times J$ with $P(\Td\times J)P^{-1}$, we write the identity
$$\Td\times J=LC^T_pL^{-1}$$
Then, by duality, we have
$$LC_p^TL^{-1}=\Ad K^{-1}C_p^TKJ'\Ad$$
where we have set the LU-decomposition $$H_n(q/p)=KJ'K^T.$$ 
Let us introduce $\bar{J'}=\Ad J'\Ad$ ; we get  :
$$\left(LC_p^TL^{-1}J\right)\times (J\bar{J'})=\Ad K^{-1}C_pKJ' \Ad$$
We remark that the matrices $LC_p^TL^{-1}J$ and $K^{-1}C_pKJ'$ are both tridiagonal and symmetric with non-singular principal diagonals, so that we necessarily have $$J\bar{J'}=\pm\Id.$$ Notice that by assumption the last coefficient of $J$ is $+1$ and that the first coefficient of $J'$ is always $+1$ (since it is the leading coefficient of $\frac{q(x)}{p(x)}$). Thus
$$J\bar{J'}=\Id.$$
Then, we may conclude by Proposition \ref{real_counting}.\par An alternative way to make use of this computation is to say that we get another proof of the equality 
 $$\PmV(\nu)=\sgn(\Bez(p,q))$$
which appears in the sequence of identities  
 $$\sgn(\Bez(p,q))=\sgn(H_n(q/p))=\sgn(J')=\sgn(J)=\PmV(\nu)=\TaQ(r,p).$$

\begin{rem}
It is possible to extend Theorem \ref{realroots_tridiag} in the case where principal diagonals of $\Td=\Td(\alpha,\beta,\beta)$ are singular. Namely, for all $k$ such that $\beta_k=0$, we  have to assume that the corresponding $k$-th entry onto the diagonal of $J$ is equal to $+1$. Then, we get that the number of real roots of $p(x)$, {\it counted with multiplicity}, is greater than $\sgn(J)$.\par
To see this, it suffices to note that the polynomial defined by $p(x)=\det(J)\det(xJ-\Td)$ factorizes through $$p(x)=\det(J_1)\det(xJ_1-\Td_{k})\times\det(J_2)\det(xJ_2-\overline{\Td}_{n-k})$$
Moreover, the matrices $\Td_{k}$ and $\overline{\Td}_{n-k}$ remain tridiagonal symmetric and $J_1$, $J_2$ remain signature matrices. If we denote by $\bigoplus$ the usual direct sum of matrices, we have $J=J_1\bigoplus J_2$ and $\Td=\Td_{k}\bigoplus\overline{\Td}_{n-k}$.\par Thus, we may proceed by induction on the degree of $p(x)$.
\end{rem}

Before stating the converse property of Theorem \ref{realroots_tridiag}, we establish a genericity lemma.

\begin{lem}\label{genericity}
Let $p(x)$ be a monic polynomial of degree $n$ with only single roots  and $q(x)=x^{n-1}+b_{1}x^{n-1}+\ldots+b_{n-1}$.
Then, the set of all $(n-1)$-tuples $(b_1,\ldots, b_{n-1})\in\R^{n-1}$ such that there is an integer $k\in\{1,\ldots, n\}$ satisfying $\det(H_k(q/p))=0$, is a proper subvariety of $\R^{n-1}$. 
\end{lem}
\begin{proof}
We only have to show that for all $k$, $\det(H_k(q/p))$, viewed as a polynomial in the variables $b_1,\ldots, b_{n-1}$, is not the zero polynomial.\par
Let $H_n(q/p)=(s_{i+j})_{0\leq i,j\leq n-1}$ where $$\frac{q(x)}{p(x)}=\sum_{j=0}^\infty\frac{s_j}{x^{j+1}}$$ and denote by $\alpha_1,\ldots,\alpha_n$ the set of all (possibly complex) roots of $p(x)$. Then, 
$$s_j=\sum_{i=1}^n\frac{q(\alpha_i)}{p'(\alpha_i)}\alpha_i^j$$
Let us introduce the real numbers defined as $$u_j=\sum_{i=1}^n\frac{\alpha_i^j}{p'(\alpha_i)}$$ 
We obviously have $u_j=0$ whenever $j\leq n-2$ and also $u_{n-1}=1$ (look at $\lim_{x\to +\infty}\frac{x^jq(x)}{p(x)}$). So that we deduce :\air
$\left\{\begin{array}{l} s_0=1\\ s_1=b_1+u_{n}\\ 
{\rm and\; more\; generally}\\
(\forall j\in\{1,\ldots,2n-2\})\; (s_j=b_j+b_{j-1}u_{n}+\ldots+b_1u_{n+j-2}+u_{n+j-1})
\end{array}\right.$\air

Then, it becomes clear that $H_{k+1}(q/p)\not \equiv 0$ for any $k$ such that $k\leq \lfloor\frac{n-1}{2}\rfloor=r$, since $s_{2k}\in\R[b_1,\ldots,b_{2k}]$ has degree $1$ in the variable $b_{2k}$ and so is the case for $H_{k+1}(q/p)$.
\par
Next, for $r<k\leq n$, we develop the determinant $H_k(q/p)$ successively according to the first columns, and we remark that its degree in the variable $b_{n-1}$ is equal to $2k-n$ (with leading coefficient equal to $-1$).  
This concludes the proof.
\end{proof}

In other words, the Lemma says that the condition
$$(\forall k\in\{1,\ldots, n\})\;(\det(H_k(q/p))=0)$$
 is {\it generic} with respect to the space of coefficients of the polynomial $q(x)$. Because of the relations between coefficients and roots, the condition is also generic with respect to the (possibly complex) roots of the polynomial $q(x)$.\air

Here is our converse statement about real roots counting :
\begin{thm}
Let $p(x)$ be a monic polynomial of degree $n$ which has exactly $s$ real roots counted with multiplicity. We can find effectively a generic family of  symmetric tridiagonal matrices $\Td$ and signature matrices $J$ with $\sgn(J)=s$, and such that $$p(x)=\det(J)\times\det(xJ-\Td).$$   
\end{thm}

\begin{proof}
If $p(x)$ has multiple roots, then we may factorize it by $\gcd(p,p')$ and use the multiplicative property of the determinant to argue by induction on the degree. Now, we assume that $p(x)$ has only simple roots.\air

We take for $q(x)$ any monic polynomials of degree $n-1$ which has exactly $s-1$ real roots interlacing those of $p(x)$. Namely, if we denote by $x_1<\ldots<x_s$ all the real roots of $p(x)$ and by $y_1<\ldots <y_{s-1}$ all the real roots of $q(x)$, we ask that $x_1<y_1<x_1<y_2<\ldots<y_{s-1}<x_s$.   
\air
Let $r(x)$ be the unique polynomial of degree $<n$ such that $r(x)\equiv\frac{q(x)}{p'(x)}\mod p(x)$ (since $p'(x)$ is invertible modulo $p(x)$).\par
From $p'r\equiv q\mod p$ and $p'(x_i)=q(x_i)$ for all real root $x_i$ of $p(x)$, we get $$\TaQ(r,p)=s=\#\{x\in\R\mid p(x)=0\}$$ 

\air
At  this point, we need that $q(x)$ satisfies another hypothesis : that is $\prs(p,q)$ shall not have any degree breakdown, or equivalently that $H(q/p)$ shall admit a $LU$-decomposition $H_n(q/p)=KJK^T$. According to Lemma \ref{genericity}, this hypothesis is generically satisfied, although it may not be always satisfied for the natural candidate $q(x)=p'(x)/n$.\par

Then, we get from Theorem \ref{duality}
$$p(x)=\det(xJ-K^{-1}C_p^TKJ)$$
where $\Td=K^{-1}C_p^TKJ$ is tridiagonal symmetric and $J$ is a signature matrix.

\air
By the proof of Proposition \ref{realroots_tridiag}, we get moreover that $$\sgn(J)=\TaQ(r,p)=\sgn(H_n(q/p)).$$
This concludes the proof since $\TaQ(r,p)=s$. 
\end{proof}

\begin{rem}
\ib
\item[(i)] The choice of such polynomials $q(x)$ with the interlacing roots property need to count and localize the real roots of $p(x)$. It can be done via Sturm sequences for instance.
\item[(ii)]
Although the polynomial $q(x)=p'(x)/n$ has not necessarily the interlacing property in general, it is the case when all the roots of $p(x)$ are real and simple. 
Moreover, in this case, the  interlacing roots condition is equivalent to the no degree breakdown condition. Indeed, $\TaQ(p'q \mod p,p)=n$ if and only if $p'(x)$ and $q(x)$ have same signs at each root of $p(x)$.
\ie
\end{rem}


\section{Some worked examples}
In order to get lighter formulas in our examples, we decide to get rid off denominators. That is why we replace signature matrices by only non-singular diagonal matrices. If one wants to deduce formulas with signature matrices, it suffices to normalize.
 
\ib
\item[1)] 
Let $p(x)=x^3+s x+t$ with $s\not =0$, and $q(x)=p'(x)=3x^2+s$. Let us introduce the discriminant of $p(x)$ as $\Delta=-4s^3-27t^2$. Consider the decomposition of the Hankel matrix

$$H(q/p)=\left(\begin{array}{ccc}
3&0&-2s\\
0&-2s&-3t\\
-2s&-3t&2s^2
\end{array}\right)
=KJK^T$$
where $$K=
\left(\begin{array}{ccc}
1&0&0\\
0&1&0\\
-\frac{2s}{3}&\frac{3t}{2s}&1\\
\end{array}\right)$$ and 
$$J=\left(\begin{array}{ccc}
3&0&0\\
0&-2s&0\\
0&0&\frac{-\Delta}{6s}\\
\end{array}\right)$$

We recover the well-known fact that $p(x)$ has three distinct real roots if and only if $s<0$ and $\Delta>0$, which obviously reduces to the single condition $\Delta>0$.\par
Then, we have the determinantal representation 

$$\Delta \times p(x)=\det(xJ-\Td)$$
where
$$\Td=\left(\begin{array}{ccc}
0&-2s&0\\
-2s&-3t&\frac{-\Delta}{6s}\\
0&\frac{-\Delta}{6s}&\frac{t\Delta}{4s^2}\\
\end{array}\right)$$

\item[2)] Consider the polynomial $p(x)=x^5-5x^3+4x$, which in fact factorizes through $p(x)=x(x-1)(x+1)(x-2)(x+2)$. Let $q(x)=p'(x)/5$. We have
$$N(q/p)=\left(
\begin{array}{ccccc}
5&0&10&0&34\\
0&10&0&34&0\\
10&0&34&0&130\\
0&34&0&130&0\\
34&0&130&0&514\\
\end{array}\right)$$

$$\Td=\left(\begin{array}{ccccc}
0&\sqrt{2}&0&0&0\\
\sqrt{2}&0&\sqrt{\frac{7}{5}}&0&0\\
0&\sqrt{\frac{7}{5}}&0&\sqrt{\frac{36}{35}}&\\
0&0&\sqrt{\frac{36}{35}}&0&\sqrt{\frac{4}{7}}\\
0&0&0&\sqrt{\frac{4}{7}}&0\\
\end{array}\right)$$

$$p(x)=\det(x\Id-\Td).$$

In order to get some parametrized identities, let us introduce the following family of polynomials $$q_a(x)=(x-a)\left(x+\frac{3}{2}\right)\left(x+\frac{1}{2}\right)\left(x-\frac{1}{2}\right).$$

We write the LU-decomposition 
$$H(q_a/p)=\left(\begin{array}{ccccc}
1&\frac{3}{2}-a&\frac{-3a}{2}+\frac{19}{4}&\frac{57}{8}-\frac{19a}{4}&\frac{-57a}{8}+\frac{79}{4}\\
&&&&\\
\frac{3}{2}-a&\frac{-3a}{2}+\frac{19}{4}&\frac{57}{8}-\frac{19a}{4}&\frac{-57a}{8}+\frac{79}{4}&\frac{237}{8}-\frac{79a}{4}\\
&&&&\\
\frac{-3a}{2}+\frac{19}{4}&\frac{57}{8}-\frac{19a}{4}&\frac{-57a}{8}+\frac{79}{4}&\frac{237}{8}-\frac{79a}{4}&\frac{-237a}{8}+\frac{319}{4}\\
&&&&\\
\frac{57}{8}-\frac{19a}{4}&\frac{-57a}{8}+\frac{79}{4}&\frac{237}{8}-\frac{79a}{4}&\frac{-237a}{8}+\frac{319}{4}&\frac{957}{8}-\frac{319a}{4}\\
&&&&\\
\frac{-57a}{8}+\frac{79}{4}&\frac{237}{8}-\frac{79a}{4}&\frac{-237a}{8}+\frac{319}{4}&\frac{957}{8}-\frac{319a}{4}&\frac{-957a}{8}+\frac{1279}{4}\\
\end{array}\right)=K_aJ_aK_a^T$$

where the associated ``signature" matrix $J_a$ is equal to 
{\footnotesize $$\left(\begin{array}{ccccc}
1&&&&\\
&-\frac{1}{2}(a+1)(2a-5)&&&\\
&&\left(\frac{15}{16}\right)\frac{(2a-1)(4a^2-a-15)}{(a+1)(2a-5)}&&\\
&&&\left(\frac{45}{128}\right)\frac{48a^4-16a^3-216a^2+58a+105}{(2a-1)(4a^2-a-15)}&\\
&&&&\left(\frac{315}{8}\right)\frac{(a+2)(a+1)a(a-1)(a-2)}{48a^4-16a^3-216a^2+58a+105}\\
\end{array}\right)$$}

The condition for $H_a(q/p)$ to be positive definite is equivalent to having only positive coefficients onto the diagonal of $J_a$.\par
First, it yields $J_a(2,2)>0$, which means that $a\in]-1,\frac{5}{2}[$. Then, we add the condition $J_a(3,3)>0$ which means that $a\in]\frac{1}{2},2,06..[$.  Then, we add the condition $J_a(4,4)>0$ which means that $a\in]0,9..,2,00..[$. And finally, we add the condition $J_a(5,5)>0$, which means that $a\in]1,2[$ and gives exactly the interlacing property for the polynomial $q_a(x)$. 

\air
For instance, with $a=\frac{3}{2}$ we get $p(x)=\det\left(x\Id-\Td_{\frac{3}{2}}\right)$ where :

$$\Td_{\frac{3}{2}}=\left(\begin{array}{ccccc}
0&\sqrt{\frac{5}{2}}&0&0&0\\
\sqrt{\frac{5}{2}}&0&\sqrt{\frac{9}{8}}&0&0\\
0&\sqrt{\frac{9}{8}}&0&\sqrt{\frac{35}{40}}&\\
0&0&\sqrt{\frac{35}{40}}&0&\sqrt{\frac{1}{2}}\\
0&0&0&\sqrt{\frac{1}{2}}&0\\
\end{array}\right)$$

\ie


\air\noindent{\small{\bf Acknowledgments.} }\par I wish to thank Marie-Francoise Roy for helpful discussions on the subject.



\end{document}